\documentclass[12pt]{article}
\usepackage[centertags]{amsmath}
\usepackage{amsfonts}
\usepackage{amssymb}
\usepackage{amsthm}
\usepackage{amsmath,mathrsfs}
\usepackage{amsmath,amscd}
\usepackage[matrix,arrow,curve,cmtip]{xy}
\title{\textbf{Goodwillie Calculus and Geometric Stacks}}
\author{Renaud Gauthier }
\theoremstyle{definition}
\newtheorem*{meta}{Theorem}

\DeclareMathOperator*{\colim}{\text{colim}}

\DeclareMathOperator*{\hocolim}{\text{hocolim}\,}
\DeclareMathOperator*{\holim}{\text{holim}}
\DeclareMathOperator*{\tildebar}{\simeq}

\newcommand{\beq}{\begin{equation}}
\newcommand{\eeq}{\end{equation}}

\newcommand{\hrarr}{\hookrightarrow}
\newcommand{\rarr}{\rightarrow}

\newcommand{\larr}{\leftarrow}
\newcommand{\Rarr}{\Rightarrow}
\newcommand{\eset}{\emptyset}
\newcommand{\Ob}{\text{Ob\,}}

\newcommand{\xrarr}{\xrightarrow}

\newcommand{\cA}{\mathcal{A}}

\newcommand{\cC}{\mathcal{C}}

\newcommand{\cD}{\mathcal{D}}
\newcommand{\cE}{\mathcal{E}}
\newcommand{\cF}{\mathcal{F}}

\newcommand{\cH}{\mathcal{H}}

\newcommand{\cM}{\mathcal{M}}

\newcommand{\cO}{\mathcal{O}}
\newcommand{\cP}{\mathcal{P}}

\newcommand{\cS}{\mathcal{S}}

\newcommand{\cU}{\mathcal{U}}

\newcommand{\cX}{\mathcal{X}}

\newcommand{\bG}{\mathbb{G}}

\newcommand{\bL}{\mathbb{L}}

\newcommand{\bR}{\mathbb{R}}
\newcommand{\bS}{\mathbb{S}}
\newcommand{\thT}{\mathbb{T}}
\newcommand{\bU}{\mathbb{U}}

\newcommand{\Hom}{\text{Hom}}
\newcommand{\Ho}{\text{Ho}\,}

\newcommand{\Map}{\text{Map}}

\newcommand{\op}{\text{op}}

\newcommand{\Set}{\text{Set}}
\newcommand{\Spec}{\text{Spec\,}}

\newcommand{\uHom}{\underline{\Hom}}

\newcommand{\AMod}{A\text{-Mod}}

\newcommand{\AffC}{\cA \text{ff}_{\cC}}
\newcommand{\AffChat}{\AffC^{\, \wedge}}
\newcommand{\AffCtildetau}{\AffC^{\,\sim,\tau}}

\newcommand{\cHom}{\cH \text{om}}

\newcommand{\Comm}{\text{Comm}}

\newcommand{\Dn}{\Delta^n}

\newcommand{\dashMod}{\text{-Mod}}

\newcommand{\Der}{\mathbb{D}\text{er}}
\newcommand{\DerA}{\mathbb{D}\text{er}_A}

\newcommand{\del}{\partial}

\newcommand{\ettop}{\acute{e}\text{t}}
\newcommand{\eps}{\epsilon}

\newcommand{\kMod}{k\text{-Mod}}
\newcommand{\kAff}{k\text{-Aff}}

\newcommand{\kDAff}{\text{k-}D^{-}\text{Aff}}

\newcommand{\kAlg}{k\text{-Alg}}

\newcommand{\LBA}{\mathbb{L}_{B/A}}
\newcommand{\LA}{\mathbb{L}_A}

\newcommand{\LB}{\mathbb{L}_B}

\newcommand{\oP}{\oplus}

\newcommand{\oT}{\otimes}
\newcommand{\oTA}{\oT_A}
\newcommand{\oTAL}{\oTA^{\mathbb{L}}}

\newcommand{\RuHom}{\bR \uHom}

\newcommand{\sPr}{\text{sPr}}

\newcommand{\Sh}{\text{Sh}}

\newcommand{\SetD}{\Set_{\Delta}}

\newcommand{\St}{\text{St}}

\newcommand{\Sp}{\text{Sp}}

\newcommand{\skMod}{\text{s}k\text{-Mod}}

\newcommand{\skAlg}{\text{s}k\text{-Alg}}

\newcommand{\uh}{\underline{h}}

\newcommand{\ucHom}{\underline{\cHom}}

\newcommand{\unpone}{\underline{n+1}}
\newcommand{\kDAfftildeet}{\kDAff^{\, \sim \, , \, \acute{e}t.}}
\newcommand{\BMod}{B \dashMod}
\newcommand{\Ker}{\text{Ker}}
\newcommand{\dashskAlg}{\text{-} \skAlg}
\newcommand{\kDAffhat}{\kDAff^{\, \wedge}}
\newcommand{\kAffhat}{\kAff^{\, \wedge}}
\newcommand{\kAfftildeet}{\kAff^{\,\sim \, , \, \acute{e}t.}}
\newcommand{\DStk}{D^{-}\Stk}
\newcommand{\Stk}{\St(k)}
\newcommand{\dashkAlg}{\text{-}\kAlg}
\newcommand{\uSpecA}{\uSpec A}
\newcommand{\uSpecB}{\uSpec B}
\newcommand{\LFx}{\bL_{F,x}}
\newcommand{\LFy}{\bL_{F,y}}
\newcommand{\uSpec}{\underline{\Spec}}
\newcommand{\RuSpec}{\bR \uSpec}
\newcommand{\TFx}{\mathbb{T}_{F,x}}

\newcommand{\bD}{\mathbb{D}}
\newcommand{\Aln}{A_{\leq n}}
\newcommand{\Alnmone}{A_{\leq n-1}}
\newcommand{\DqdA}{\Delta^q \downarrow A}
\newcommand{\DdA}{\Delta \downarrow A}
\newcommand{\skn}{\text{sk}_n}
\newcommand{\LFGx}{\bL_{F/G,x}}
\newcommand{\LFGy}{\bL_{F/G,y}}
\newcommand{\dash}{\text{-}}

\begin{document}
\maketitle

\begin{abstract}
	We show Goodwillie's calculus of functors and $n$-geometric $D^{-}$-stacks share similar features by starting to focus on the convergence of Taylor towers for homotopy functors and the fact that $\bR F(A) \cong \holim \bR F(\Aln)$ for geometric stacks, where $\{\Aln \}$ provides a Postnikov tower of some given $A \in \skAlg$. From there we show we have parallel results, such as similar homotopy fibers of connecting maps in towers, as well as polynomial approximations, pointwise approximations and reconstruction theorems for towers, all of which point to the presence of a common underlying theory. 
\end{abstract}

\newpage

\section{Introduction}
The present work provides an illustration of a (putative) higher categorical generalization of the bridge technique of \cite{OC} relating Morita equivalent theories $\thT$ and $\thT'$ via the equivalence of their respective classifying topos $\cE_{\thT}$ and $\cE_{\thT'}$ as in:
\beq
\xymatrix{
	& \cE_{\thT} \simeq \cE_{\thT'} \ar@{.>}[dr] \\
	\thT \ar@{.>}[ur] && \thT'
} \nonumber
\eeq
We regard theories as being inert objects in contrast to their classifying topoi which we regard as being dynamical representations thereof. Indeed, for a geometric theory $\thT$, its classifying topos $\cE_{\thT}$ is  $\Sh(\cC_{\thT}, J_{\thT})$, the Grothendieck topos of sheaves on its syntactic category $\cC_{\thT}$. This provides a dynamic manifestation of $\thT$ insofar as it is functorial. On can generalize this to the higher categorical setting by considering functors valued in $\SetD$, $\Sp$, the category of spectra, or $\cS$, the category of spaces. In particular, in the present paper we consider the category of homotopy functors $\cF(\cC, \cD)$ from $\cC$ to $\cD$, both being categories of spaces with a notion of weak equivalence, and $\kDAfftildeet$, the category of $D^{-}$-stacks on $\skAlg$ for the \'etale topology, both of which are higher categories. Goodwillie calculus has homotopy functors $F: \cC \rarr \cD$ as objects of study, and such a ``theory", in the pedestrian sense, is built upon a theory that we denote by $\bG$. In the same manner, derived stacks are functors $F: \skAlg \rarr \SetD$ satisfying certain conditions, objects of $\kDAfftildeet$, whose formalism is developed following a theory that we will denote by $\bS$. Whether it be $\bG$ or $\bS$, those are abstractions of the defining features of the respective higher categorical constructs they are meant to govern. \\

We argue different formalisms, such as Goodwillie calculus and the theory of $n$-geometric $D^{-}$-stacks, without necessarily being equivalent, at least share some attributes, which we refer to as aspects, collectively denoted by $a$, and we write 
\beq
\cF(\cC, \cD) \tildebar^a \kDAfftildeet \nonumber
\eeq
to indicate that Goodwillie calculus in $\cF(\cC, \cD)$ and $n$-geometric $D^{-}$-stacks in $\kDAfftildeet$ share some common features, namely those in $a$. The object of the present work is to provide evidence for the existence of such features. Ultimately, our aim is to show that it follows there is an underlying theory $\thT_0$, from which $\bG$ and $\bS$ both originate. In the bridge technique formalism, this may be presented as follows:
\beq
\xymatrix{
	\cF(\cC, \cD) \ar[rr]^{\tildebar^a} &&\kDAfftildeet \\
	\bG_a \ar@{.>}[u] && \bS_a  \ar@{.>}[u] \\
	& \thT_0 \ar[ul] \ar[ur]
} \nonumber
\eeq
where $\thT_a$ is that part of $\thT$ necessary to develop those aspect $a$ of a formalism $\cE_{\thT}$. In the above diagram $\bG_a$ and $\bS_a$ are seen as expansions of $\thT_0$. Observe that a theory $\thT$ can be seen as an expansion of $\thT_a$, thus the bottom part of the above diagram can be completed as such:
\beq
\xymatrix{
	\bG && \bS \\
	\bG_a \ar[u] && \bS_a \ar[u] \\
	& \thT_0 \ar[ul] \ar@{.>}[uul] \ar@{.>}[uur] \ar[ur]
} \nonumber
\eeq
a maximal case corresponding to considering all aspects of both formalisms, in which case one recovers a higher categorical generalization of the bridge technique per se.\\

Though written with the aim to extract a theory $\thT_0$ from $\bG$ and $\bS$, those latter theories are expansions of, it is not until the very end of this paper that we will discuss theories proper. Throughout the paper, when we use the word ``theory", we really mean formalism.\\

In a first time, we draw comparisons between certain results pertaining to the calculus of functors of Goodwillie (\cite{Go1}, \cite{Go2}, \cite{Go3}) and some aspects of the theory of $n$-geometric $D^{-}$-stacks (\cite{TV}, \cite{TV4}), rather than being a direct application of the former theory to the latter. Thus our aim is to show that those two different theories share some common features, something we would like to abstract further. We focus in particular on the Taylor towers for homotopy functors, on their built-in polynomials and derivatives thereof on the one hand, and on the convergence of towers built from truncations $\bR F(\Aln)$ for $A \in \skAlg$, $F$ a $n$-geometric $D^{-}$-stack on the other. After having drawn some comparisons, we study a few parallel results, namely the similitudes between the construction of polynomial approximations, point-wise equivalences, and a reconstruction theorem for towers. In the last part of the present work we go a little deeper into each theory, which allows us to formally abstract their common features into an overarching theory.\\

Goodwillie observed that to a homotopy functor $F: \cC \rarr \cD$ between categories of spaces, that is a functor that preserves equivalences, one can associate a Taylor tower $\cdots \rarr P_n F(X) \rarr \cdots \rarr P_1 F(X) \rarr P_0 F(X)$ for $X \in \cC$, and that under a certain condition called analyticity, this Taylor tower converges to $F(X)$. In Derived Algebraic Geometry, one starts with $\skMod$ as base symmetric monoidal model category, with $\skAlg$ as model category of commutative monoids. One can associate to any $A \in \skAlg$ a Postnikov tower $ \cdots \rarr \Aln \rarr \Alnmone \rarr \cdots \rarr A_0$, and if $F$ is a $n$-geometric $D^{-}$-stack, Theorem C.0.9 of \cite{TV4} shows one has an isomorphism $\bR F(A) \rarr \holim_n \bR F(\Aln)$, a fact that spurred the present study. It turns out both theories, though very different, share many common traits (deformation theoretic in nature) as mentioned above, thereby pointing to a common underlying theory. We obtain the following metatheorem:
\begin{meta}
	The theory $\bG$ of Goodwillie calculus and the theory $\bS$ of geometric stacks in $\kDAfftildeet$ are both expansions of a common theory $\thT_0$ of functors $F: \cC \rarr \cD$, $\cC$ a category with a notion of neighborhood, $\cD$ a category with weak equivalences, such that to objects $X$ of $\cC$ one can associate deformation theoretic towers $\cdots X_n \rarr X_{n-1} \rarr \cdots $ for which $F(X) \simeq \holim F(X_n)$ in $\Ho(\cD)$, provided one puts a deformation theoretic constraint on $F$ and $X$ lies in an appropriate neighborhood. 
\end{meta}

Goodwillie calculus was initiated in the seminal papers \cite{Go1}, \cite{Go2} and \cite{Go3} and pursued further in a host of later publications some of which the reader can find in the reference section. One paper in particular we use is \cite{AC}. For stacks we essentially use the foundational papers \cite{TV} and \cite{TV4}.\\
 
In section \ref{Goodwillie Calculus} we give the main results of Goodwillie Calculus for our purposes, and in section \ref{Geometric stacks} we formally define geometric stacks and present a representation theorem for such objects that is central to establishing a bridge with the calculus of functors. It is in sections \ref{The bridge} and \ref{Parallel results} that we present shared features of both theories; convergence of towers, analyticity and obstruction theory, similarities in polynomial approximations, pointwise approximations and reconstruction theorems. In section \ref{Taking stock} we abstract both theories further which allows to extract a common underlying theory $\thT_0$.

\section{Goodwillie Calculus} \label{Goodwillie Calculus}
We will start our study of Goodwillie Calculus by one of the main results of \cite{Go3}, namely that certain types of functors can have approximations by polynomial like approximations displayed in a Taylor-like sequence, and a sub-type of such functors have those Taylor towers converging to them.\\

To be precise, Goodwillie considers functors $F: \cC \rarr \cD$ between categories of spaces (based or unbased) and sometimes $\cD$ is even taken to be a category of spectra. In particular, one considers functors that preserve weak equivalences. Such functors are referred to as \textbf{homotopy functors}.\\

To have a Taylor expansion presupposes that a point be fixed, and the expansion is relative to some variable in a neighborhood of such a fixed point. This is implemented in the present situation by working with slice categories; if $Y$ is a fixed space, we therefore look at $\cS_{/Y}$, the category of based spaces over $Y$, or $\cU_{/Y}$, the category of unbased spaces over $Y$, which we generically denote by $\cC_Y$. The notion of neighborhood is then implemented by considering the connectedness of maps $X \rarr Y$ for objects $X \in \cC_Y$.\\

One can associate to any homotopy functor $F: \cC_Y \rarr \cD$ a \textbf{Taylor tower} of homotopy functors $P_n F$: \footnote{Thm. 1.13 \cite{Go3}}
\beq
\xymatrix{
	&F \ar[dl] \ar[d]_{p_nF} \ar[dr] \ar[drrr] \ar[drrrr]^{p_0F}\\
\cdots \ar[r] & P_n F \ar[r]_-{q_n F}  &P_{n-1} F \ar[r] &\cdots \ar[r] 
	& P_1 F \ar[r]_-{q_1 F} & P_0 F
} \nonumber
\eeq
The different constituents in this tower are constructed as follows. Let $F: \cC_Y \rarr \cD$ be a homotopy functor, where $\cC$ and $\cD$ are the category of based or unbased spaces, and $\cD$ could also be the category Sp of spectra. For $X \in \cC_Y$ and $U$ any space, one defines the \textbf{fiberwise join over $Y$} \cite{Go2} by:
\beq
X *_Y U = \hocolim(X \larr X \times U \rarr Y \times U) \nonumber
\eeq
Let $\unpone = \{1, \cdots, n \}$, and denote by $\cP(\unpone)$ the poset of subsets of $\unpone$. Let $\cP_0(\unpone) \subset \cP(\unpone)$ be the poset of non-empty subsets of $\unpone$. Define:
\beq
T_nF(X) = \holim_{U \in \cP_0(\unpone)} F(X *_Y U) \nonumber
\eeq
This provides a homotopy functor $T_nF: \cC_Y \rarr \cD$, along with a natural transformation $F \xrarr{t_n F} T_n F$ which pointwise reads:
\beq
F(X) = F(X *_Y \eset) \xrarr{t_n F(X)} \holim_{U \in \cP_0(\unpone)}F(X *_Y U) = T_n F(X) \nonumber
\eeq
We denote by $P_n F(X)$ the sequential homotopy colimit of:
\beq
F(X) \xrarr{t_n F(X)} T_n F(X) \xrarr{t_n(T_n F)(X)} T_n^2 F(X) \rarr \cdots \nonumber
\eeq
This provides yet another homotopy functor $P_n F: \cC_Y \rarr \cD$, along with a natural map $F \xrarr{p_n F} P_n F$. By construction we also have maps $q_{n,k}: T_n ^k F \rarr T_{n-1}^k F$. To see that this is the case, observe:
\beq
T_n^k F(X) \cong \holim_{(U_1, \cdots, U_k) \in \cP_0(\unpone)^k} F(X *_Y (U_1 * \cdots * U_k)) \nonumber
\eeq
Since we have an inclusion $\cP_0(\underline{n})^k \hrarr \cP_0(\unpone)^k$, it follows we have an induced map $T_n^k F \xrarr{q_{n,k}} T_{n-1}^k F$ as desired, hence an induced map $q_n F: P_n F \rarr P_{n-1} F$ of horizontal homotopy colimits. Along with the maps $F \xrarr{p_nF} P_nF$, this provides us with the desired Taylor tower.\\

We now discuss one important characteristic of the functor $P_nF$. In order to do so we have to introduce \textbf{cubical diagrams} \cite{Go2}, functors $\cX: \cP(S) \rarr \cC$, $S$ a finite set, $\cP(S)$ the poset of all subsets of $S$, $\cC$ one of our categories of spaces. For instance if $S = \{1,2\}$, we have $\cP(S) = \{ \eset, \{1 \}, \{2\}, \{1,2\} \}$, so that the diagram $\cX$ looks like:
\beq
\xymatrix{
\cX( \eset) = \lim \cX \ar[d] \ar[r] & \cX(2) \ar[d] \\
\cX(1) \ar[r] & \cX(1,2)
} \nonumber
\eeq
whence the name cubical. If $S$ has cardinality $n$, $\cX$ is referred to as a $n$-cube. Denote by $\cP_1(S) = \{T \, | \, T \subsetneq S \}$. We say $\cX$ is \textbf{co-cartesian} \footnote{Def. 1.4 \cite{Go2}} if the map:
\beq
\hocolim(\cX|_{\cP_1(S)}) \rarr \colim(\cX) \nonumber
\eeq
is a weak equivalence. We say $\cX$ is \textbf{strongly co-cartesian} \footnote{Def. 2.1 \cite{Go2}} if each of its faces of dimension greater than 2 is co-cartesian. In contrast to being co-cartesian, a cubical diagram $\cX$ is said to be \textbf{cartesian} \footnote{Def. 1.3 \cite{Go2}} if the map:
\beq
\cX(\eset) = \lim \cX \rarr \holim \cX|_{\cP_0(S)} \nonumber
\eeq
is a weak equivalence. We can now characterize our functors $P_n F$: a homotopy functor $F: \cC \rarr \cD$ is said to be \textbf{$n$-excisive} \footnote{Def. 3.1 \cite{Go2}} if for every strongly co-cartesian $(n+1)$-cubical diagram $\cX$, the composite cubical diagram $F \circ \cX: \cP(S) \rarr \cD$ is cartesian. By Theorem 1.8 of \cite{Go3}, $P_nF$ is $n$-excisive. A weaker condition referred to as stable $n$-th order excision (see below) is satisfied in particular by analytic functors. By Proposition 1.5 of \cite{Go3} $F$ $n$-stably excisive agrees with $P_nF$ (which is $n$-th order excisive) to order $n$ in a certain sense. Essentially the discrepancy between the two notions morally vanishes as $n$ tends to $\infty$: if $F$ is $\rho$-analytic (in particular $n$-stably excisive), for any $X \in \cC_Y$, such that $X \rarr Y$ is $(\rho +1)$-connected, i.e. $X$ is close enough to $Y$, then $F(X) \simeq \holim P_nF(X)$ \footnote{Thm. 1.13 \cite{Go3}}.\\

The interpretation of $P_n F$ as a Taylor polynomial however is not immediate. To see this clearly one has to consider the homotopy fiber of the map $P_n F(X) \rarr P_{n-1}F(X)$ for $X \in \cC$, which is denoted $D_n F(X)$ in \cite{Go3}. This should correspond to $\frac{f^{(n)}(a)}{n!} (x-a)^n$ in Calculus. To see that we have something similar in form, it is convenient to work with spectra: $\cD = \Sp$. In this case:
\beq
D_n F(X) = (C_n \wedge (X^{\wedge n}))_{h\Sigma_n} \nonumber
\eeq
where $C_n$ is a spectrum with a $\Sigma_n$-action, and $h \Sigma_n$ stands for homotopy orbit spectrum.\\

We now focus on the conditions under which we have convergence of the Taylor tower of a homotopy functor $F$. As mentioned above, by Theorem 1.13 of \cite{Go3}, if $F$ is $\rho$-analytic and $X \rarr Y$ is $(\rho + 1)$-connected, then $F(X) \simeq \holim P_nF(X)$. Being analytic is a variant on the $n$-excisiveness condition. To introduce that concept, we need a few definitions. First, a map of spaces is $k$-connected if all its homotopy fibers are $(k-1)$-connected. Recall that above we defined a cubical diagram $\cX$ to be cartesian if the map:
\beq
\cX(\eset) \rarr \holim \cX|_{\cP_0(S)} \nonumber
\eeq
is a weak homotopy equivalence. We now say $\cX$ is \textbf{$k$-cartesian} \footnote{Def. 1.3 \cite{Go2}} if this map is $k$-connected. All of this being defined, we can define a weaker form of $n$-excisiveness: $F$ as a functor is \textbf{stably $n$-excisive} \footnote{Def. 4.1 \cite{Go2}} if it satisfies the following condition $E(n,c,K)$: for any strongly co-cartesian $(n+1)$-cubical diagram $\cX: \cP(S) \rarr \cC$ such that $\forall s \in S$, $\cX(\eset) \rarr \cX(s)$ is $k_s$-connected for $k_s \geq K$, the diagram $F \circ \cX$ is $(-c + \sum k_s)$-cartesian. Finally $F$ is said to be \textbf{$\rho$-analytic} \footnote{Def. 4.2 \cite{Go2}} if it satisfies $E(n,n \rho - q, \rho + 1)$ for some $q$, and this for all $n \geq 1$.

\section{Geometric stacks} \label{Geometric stacks}

\subsection{General context}
We briefly go over the construction of stacks. For a complete coverage, see \cite{TV4} and \cite{TV}. Given a homotopical algebraic context $(\cC, \cC_0, \cA)$ (to be discussed below), let $\Comm(\cC)$ be the model category of commutative monoid objects of $\cC$ with $\AffC$ as opposite category, which we refer to as the model category of affine schemes, or affine stacks. Let $\AffChat = L_W \sPr(\AffC)$ be the model category of prestacks on $\AffC$,  left Bousfield localization of the model category of simplicial functors on $\AffC$ with respect to weak equivalences in $\AffC$. After having fixed a model topology $\tau$ on this latter model category, one can then define the model category of stacks $\AffCtildetau$ on $\AffC$ as the left Bousfield localization of $\AffChat$ with respect to homotopy $\tau$-hypercovers.\\

For the sake of defining geometric stacks, one needs to introduce several assumptions, collected together into what is called a \textbf{Homotopical Algebraic context} (referred to as HA), and a \textbf{Homotopical Algebraic Geometry} context (referred to as HAG). For a full list of assumptions, the reader is referred to the main reference \cite{TV4} where those are introduced. We will just present those that are most useful for our purposes. A HA context is a triple $(\cC, \cC_0, \cA)$ consisting of a symmetric monoidal model category $(\cC, \otimes, 1)$. If we fix universes $\mathbb{U} \in \mathbb{V}$, we suppose $\cC$ is $\mathbb{V}$-small, and $\mathbb{U}$-combinatorial. It is also assumed proper, pointed, and $\Ho(\cC)$ is an additive category. For $A \in \Comm(\cC)$, $A \dashMod$ is a $\mathbb{U}$-combinatorial model category, and it is proper. One can define a tensor product on $A \dashMod$, for which it is a symmetric monoidal model category. $\cC_0$ is a full subcategory of $\cC$ that one needs in obstruction theory, it is necessary to define various kinds of geometric stacks. $\cA$ is a full subcategory of $\Comm(\cC)$ satisfying a technical condition relative to $\cC_0$.\\

A HAG context consists of a HA context, together with a choice of a model topology $\tau$ on $\AffC$, and a class of morphisms $P$ in $\AffC$ which is stable by equivalences. There are further assumptions on $\tau$ and $P$ so that $(\cC, \cC_0, \cA, \tau, P)$ becomes a HAG context, some more technical than others. One can mention that covering families use morphisms in $P$, and morphisms in this class are stable by equivalences, homotopy pullbacks and compositions. The topology is necessary to define stacks, and the choice of $P$ will dictate what kind of geometric stack we get, assuming morphisms in $P$ are compatible with the topology.\\

We suppose fixed a HAG context $(\cC, \cC_0, \cA, \tau, P)$. We can now define the notion of \textbf{$n$-geometric stack} \footnote{Def. 1.3.3.1 \cite{TV4}} using induction: a stack $F \in \AffCtildetau$ is \textbf{$n$-geometric} if the diagonal morphism $F \rarr F \times^h F$ is $(n-1)$-representable, and if $F$ admits an $n$-atlas. One defines a \textbf{$n$-atlas} for $F$ to be a $\mathbb{U}$-small family of morphisms $\{U_i \rarr F \}_{i \in I}$ for which each $U_i$ is representable, $\coprod U_i \rarr F$ is an epi, and each morphism in this family is in $(n-1) \dash P$. One says a morphism of stacks $F \rarr G$ \textbf{is in $n \dash P$} if it is $n$-representable, and for any morphism of stack $X \rarr G$, there is a $n$-atlas $\{U_i\}$ for $F \times^h_G X$ such that each composite $U_i \rarr X$ is in $P$. Finally a morphism of stacks $F \rarr G$ is \textbf{$n$-representable} if for any morphism of stack $X \rarr G$, $F \times^h_G X$ is $n$-geometric. To initialize the induction, a representable stack is said to be $(-1)$-geometric, a morphism of stacks $F \rarr G$ is $(-1)$-representable if for any morphism $X \rarr G$ from a representable stack $X$, $F \times^h_G X$ is representable as well, i.e. it is locally representable, and $F \rarr G$ is said to be in $(-1) \dash P$ if in addition the morphism $F \times^h_G X \rarr X$ is in $P$, that is $F \rarr G$ is locally in $P$. Observe that it follows from Proposition 1.3.3.3 of \cite{TV4} that a $(n-1)$-geometric stack is $n$-geometric.\\

We are interested in particular in the setting where $\cC= \skMod = \cC_0$, $\cA = \skAlg$, $\tau = \ettop$, and $P$ is the class of smooth morphisms which we define below. This provides us with a HAG context, which allows us to talk about \textbf{$n$-geometric $D^{-}$-stacks}. We denote $\Comm(\skMod)^{\op} = \kDAff$, with the homotopy category of $D^{-}$-stacks being denoted $\DStk = \Ho(\kDAfftildeet)$.\\

We say a morphism $f:A \rarr B$ of commutative algebras is \textbf{smooth} \footnote{Def. 1.2.7.2 \cite{TV4}} if it is finitely presented in $\skAlg = \Comm(\skMod)$ and it is \textbf{formally smooth} \footnote{Def. 1.2.7.1 \cite{TV4}}, that is $\LBA \in \BMod$ is projective, and there is a retraction of the morphism $\LA \oTAL B \rarr \LB$ in $\Ho(\BMod)$. Recall that for $A \in \skAlg$, to say $B \in \AMod$ is \textbf{projective} \footnote{Def. 1.2.4.1 \cite{TV4}}means it is a retract of $\coprod_{S}^{\bL} A$ for some small $\bU$-set $S$. The $B$-module $\LBA$ is the \textbf{cotangent complex of $B$ over $A$} \footnote{Def. 1.2.1.5 \cite{TV4}}, object of $\Ho(\BMod)$, constructed as:
\beq
\LBA = \bL Q \bR I( B \oTAL B) \nonumber
\eeq
with $Q: B - \skAlg_{\text{nu}} \rarr \BMod$ is the morphism from non-unital $B$-algebras defined by the pushout:
\beq
\xymatrix{
	C \oT_B C \ar[d] \ar[r]^{\mu} & C \ar[d] \\
	\text{*} \ar[r] & Q(C)
} \nonumber
\eeq
and $I : B -\skAlg_{/B} \rarr B-\skAlg_{\text{nu}}$ takes a diagram $B \xrarr{f} C \xrarr{g} B$ to $\Ker(g)$. \\

Note that if $A = 1$, we let $\bL_{B/1} = \LB$, and we refer to it as the cotangent complex of $B$ \footnote{Def. 1.2.1.5 \cite{TV4}}.\\

To justify the existence of $\LBA$, we have to introduce another definition; if $A \rarr B$ is a morphism in $\skAlg$, and $M \in \BMod$, one defines the \textbf{simplicial set of derived $A$-derivations}: \footnote{Def. 1.2.1.1 \cite{TV4}} $B \rarr M$ to be defined by:
\beq
\DerA(B,M) = \Map_{A\dashskAlg_{/B}}(B, B \oP M) \nonumber
\eeq
an object of $\Ho(\SetD)$, where $B \oP M \in \skAlg$ has $A \coprod M$ as underlying object, with product as defined in \cite{TV4}. Now we have the following result from \cite{TV4} which guarantees the existence of a cotangent complex: \footnote{Prop. 1.2.1.2 \cite{TV4}} given a morphism $A \rarr B$ in $\skAlg$, there is some $\LBA \in \BMod$ and some $d \in \pi_0(\DerA(B,\LBA))$ such that $\forall \, M \in \BMod$, the natural map $d^*: \Map_{\BMod}(\LBA,M) \rarr \DerA(B,M)$ induced by $d$ is an isomorphism in $\Ho(\SetD)$.\\

An alternate characterization of smooth morphisms in $\skAlg$ is provided by Theorem 2.2.2.6 of \cite{TV4}: a morphism is smooth if and only if it is \textbf{strongly smooth}\footnote{Def. 2.2.2.3 \cite{TV4}}, that is, it is \textbf{strong} (for all $k$, $\pi_k(A) \oT_{\pi_0(A)} \pi_0(B) \rarr \pi_k(B)$ is an isomorphism) and $\Spec \pi_0 B \rarr \Spec \pi_0 A$ is smooth.\\

\subsection{Representability Theorem}
Central to our transition from Goodwillie calculus to its essential realization within the realm of derived stacks is the following representability Theorem, a pared down version of the original result of Lurie (\cite{Lu}) adapted to the setting of geometric stacks by Toen and Vezzosi (\cite{TV4}): a $D^{-}$-stack $F$ is $n$-geometric if and only its truncation $t_0(F)$ is an Artin $(n+1)$-stack, $F$ has an obstruction theory, and $\forall \, A \in \skAlg$, we have a natural isomorphism $\bR F(A) \rarr \holim_k \bR F(A_{\leq k})$ in $\Ho(\SetD)$. We discuss these conditions in turn.\\

\subsubsection{Truncation statement}
For the \textbf{truncation functor} \footnote{Def. 2.2.4.3 \cite{TV4}}, we have a natural inclusion functor $i: \kAff \hrarr \kDAff$, inducing a right Quillen functor $i^*: \kDAffhat \rarr \kAffhat$, with a right derived functor $\bR i^*: \Ho(\kDAffhat) \rarr \Ho(\kAffhat)$ 
that preserves the sub-category of stacks, hence provides a right Quillen restriction functor $i^*: \kDAfftildeet \rarr \kAfftildeet$, whose right derived functor $t_0 = \bR i^*: \DStk \rarr \Stk$ is our desired truncation functor. Here $\Stk$ is the category of stacks obtained by starting with $\kMod$ instead of $\skMod$.\\

An \textbf{Artin $n$-stack} \footnote{Def. 2.1.1.4 \cite{TV4}} is an object of $\Stk$ that is $n$-truncated and $p$-geometric for some integer $p$. Here the HAG context is based on $\cC = \kMod = \cC_0$, $\cA = \kAlg$, $\tau = \ettop$ and $P$ is the class of smooth morphisms, where we say $A \rarr B$ in $\kAlg$ is \textbf{smooth} if $B \in A \dashkAlg$ is finitely presented, $\Omega_{B/A}$ is projective in $\BMod$, and $\pi_k(\LBA) = 0$ for $k \geq 1$. To say that $t_0(F)$ is an Artin $(n+1)$-stack is due to the fact that if $F$ is a $n$-geometric stack in $\St(k)$, it is $(n+1)$-truncated.\\

\subsubsection{Obstruction theory statement}
The notion of obstruction theory is slightly more involved. For a pointed model category $\cC$, we have a suspension functor $ \Sigma: \Ho(\cC) \rarr \Ho(\cC)$, $x \mapsto * \coprod_x^{\bL} * $. Denote by $\cC_1 = \Sigma \cC_0$ the full subcategory of $\cC$ spanned by objects equivalent to the suspension of some object of $\cC_0$. The obstruction for $n$-geometric $D^{-}$-stacks is relative to $\skMod_1$.\\

A stack has an \textbf{obstruction theory} \footnote{Def. 1.4.2.1. \cite{TV4}}relative to a HA context $(\cC, \cC_0, \cA)$ if it is infinitesimally cartesian and if it has a global cotangent complex.\\

A stack $F$ is said to be \textbf{infinitesimally cartesian} \footnote{Def. 1.4.2.1 \cite{TV4}}relative to a given HA context $(\cC, \cC_0, \cA)$, if $\forall \, A \in \cA$, $M \in \AMod_1 = \Sigma \AMod_0$, $d \in \pi_0(\DerA(A,M))$ with a corresponding morphism $A \rarr A \oP M$, still denoted by $d$, the square:
\beq
\xymatrix{
	\bR F(A \oP_d \Omega M) \ar[d] \ar[r] & \bR F(A) \ar[d]^d \\
	\bR F(A) \ar[r]_-s & \bR F(A \oP M)
} \nonumber
\eeq
is homotopy cartesian. Recall what this means: from \cite{GoJa} for instance, if $\cM$ is a right proper model category, $A,B,C,D \in \Ob \cM$, a square:
\beq
\xymatrix{
	A \ar[d] \ar[r] & C \ar[d] \\
	B \ar[r] & D
} \nonumber
\eeq
is said to be homotopy cartesian, if for any factorization of $C \rarr D$ into a trivial cofibration $C \rarr E(C)$ followed by a fibration $E(C) \rarr D$, then the induced map $A \xrarr{i_*} B \times_D E(C)$ is a weak equivalence. Note also that since we assume $\cC$ to be pointed, it has a suspension functor as mentioned above, with the loop functor $\Omega: \Ho (\cC) \rarr \Ho(\cC)$, $x \mapsto * \times^h_x *$ as right adjoint. Then $A \oP_d \Omega M$ denotes the square zero extension of $A$ by $\Omega M$, defined by the homotopy cartesian square:
\beq
\xymatrix{
A \oP_d \Omega M \ar[d] \ar[r] & A \ar[d]^d \\
A \ar[r]_-s &A \oP M
} \nonumber
\eeq
where $A \in \Comm(\cC)$, $M \in \AMod$, $d: A \rarr A \oP M$ a derivation, $s: A \rarr A \oP M$ the trivial derivation. In other terms $F$ infinitesimally cartesian means that it preserves this homotopy square. To explain the notation, denote $A \oP_d \Omega M$ by $B$ for the moment. The morphism $B \rarr A $ has a fiber isomorphic to $\Omega M = * \times^h_M *$, with the zero multiplication, and that is why $B$ is called a zero extension of $A$ by $\Omega M$. Recall the definition of such objects; for a commutative ring $R$, a square zero extension of $R$ is given by a surjective morphism of commutative rings $\phi: R' \rarr R$ such that $(\Ker \phi)^2 = 0$, which is what we have here.\\

\newpage

A stack $F$ is said to have a \textbf{global cotangent complex} \footnote{Def. 1.4.1.7 \cite{TV4}} relative to a given HA context $(\cC, \cC_0,\cA)$ if it has a local cotangent complex that behaves well under base change. More precisely, $\forall \, A \in \cA$, $x: \bR \uSpecA \rarr F$, $F$ has a (local) cotangent complex $\LFx$ at $x$, such that for any morphism $A \rarr B$ in $\cA$, for any morphism:
\beq
\xymatrix{
	\bR \uSpecB \ar[dr]_y \ar[rr]^u && \bR \uSpecA \ar[dl]^x \\
	&F
} \nonumber
\eeq
we have an induced isomorphism $u^*: \LFx \oTAL B \rarr \LFy$ in $\Ho(\Sp(\BMod))$. In this definition, one uses the Yoneda embedding $\uh: \AffC \rarr \AffChat$, with a total right derived functor $\bR \uh: \Ho(\AffC) \rarr \Ho(\AffChat)$, which turns out to be a fully faithful functor into the category of stacks $\St(\cC,\tau)$. This follows from our assumption on the model topology $\tau$ to have a HAG context, which implies in particular that it be sub-canonical. It follows that we have $\bR \uh: \Ho(\AffC) \rarr \St(\cC,\tau)$, and one denotes $\bR \uh_{\Spec A}$ by $\bR \uSpecA$. Regarding the (local) cotangent complex, it is defined as follows. If $F$ is a stack, $A \in \cA$, $x: \bR \uSpecA \rarr F$, we say $F$ has a \textbf{cotangent complex} \footnote{Def. 1.4.1.5 \cite{TV4}} at $x$ if there exists some $(-n)$-connective stable $A$-module $\LFx \in \Ho(\Sp(\AMod))$ and an isomorphism $\Der_F(X,-) \cong \bR \uh_s^{\LFx}$ in $\Ho((\AMod_0^{\op})^{\wedge})$. In such a definition, the derivations functor $\Der_F(X,-)$ is defined like so: given $x : X = \bR \uSpecA \rarr F$, denote $X[M] = \bR \uSpec(A \oP M)$. One can then define the \textbf{simplicial set of derivations from $F$ to $M$ at the point $x$} \footnote{Def. 1.4.1.4 \cite{TV4}} by:
\beq
\Der_F(X,M) = \Map_{X/\AffCtildetau}(X[M],F) \nonumber
\eeq
Regarding stable modules, recall that $\cC$ is pointed, so we have suspension and loop functors, as mentioned above. Fix a cofibrant model $S^1_{\cC}$ for the object $\Sigma(1) \in \Ho(\cC)$, and for $A \in \Comm(\cC)$, define $S^1_A = S^1_{\cC} \oT A \in \AMod$. Then one defines the \textbf{model category of stable $A$-modules} \footnote{Def. 1.2.11.1 \cite{TV4}}to be the model category of spectra $\Sp(\AMod)$ with the connecting map $S^1_A \oTA - $, whose objects will simply be referred to as stable modules. \\

From there one defines a functor:
\begin{align}
	\uh_s^{-}: \Sp(\AMod)^{\op} & \rarr (\AMod_0^{\op})^{\wedge} \nonumber\\
	M_* & \mapsto \uh_s^{M_*} \nonumber
\end{align}
where:
\begin{align}
	\uh_s^{M_*}: \AMod_0  \rarr & \SetD \nonumber \\
	N  \mapsto & \Hom(M_*, \Gamma_*(S_A(N))) \nonumber \\ 
	&= \Map_{\Sp(\AMod)}(M_*, S_A(N)) \nonumber
\end{align}
with $\Gamma_*$ the simplicial resolution functor on $\Sp(\AMod)$, and $S_A: \AMod \rarr \Sp(\AMod)$ has $n$-th space defined by $S_A(N)_n = (S^1_A)^{\oTA n} \oTA N$. By Proposition 1.2.11.3 of \cite{TV4}, we have a right derived functor $\bR \uh_s^{-}: \Ho(\Sp(\AMod))^{\op} \rarr \Ho((\AMod_0^{\op})^{\wedge})$. \textbf{Connectivity} is defined by induction: $M_* \in \Ho(\Sp(\AMod))$ is said to be $0$-connective if it is isomorphic to some object of $S_A(\Ho(\AMod_0))$. Then $M_*$ is said to be $(-n)$-connective if it is isomorphic to $\Omega (M_*')$ for some $-(n-1)$-connective $A$-module $M'_*$.\\

To summarize that part of having an obstruction that involves having a cotangent complex, for a stack to have a local cotangent complex means that functorially:
\beq
\xymatrix{
\Der_F(X,-) \ar@{=}[d] \ar[r]^{\cong} &\bR \uh_s^{\LFx} \ar@{=}[d] \\
\Map_{_X/\AffCtildetau}(\bR \uSpec(A \oP -),F) & \bR \Hom(\LFx,\Gamma_*(S_A(-)))
} \nonumber
\eeq
To see that the cotangent complex does exhibit a notion of cotangency, we have to relate this notion to that of a \textbf{tangent complex} introduced in the same definition 1.4.1.5 of \cite{TV4} that the cotangent complex was introduced. By definition, the tangent complex of $F$ at the same point $x$ is defined by $\TFx = \RuHom(\LFx,A) \in \Ho(\Sp(\AMod))$, where we used the fact that $\Ho(\Sp(\AMod))$ is a closed symmetric monoidal category, hence has a stable $A$-module of morphisms $\RuHom(-,-)$. The tangent complex can be related to the \textbf{tangent stack $TF$} \footnote{Def. 1.4.1.2 \cite{TV4}} of a stack $F$. The latter is defined as $TF = \bR_{\tau} \ucHom (\bD_{\eps},F)$, where $\bR_{\tau} \ucHom$ is the internal Hom of $\Ho(\AffCtildetau)$, and $\bD_{\eps} = \RuSpec(1[\eps])$ is the infinitesimal disk, with $1 \in \Comm(\cC)$ having a trivial square zero extension $1 \oP 1 = 1[\eps]$. In this manner, $TF = \bR_{\tau} \ucHom(\bR \uSpec( 1 \oP 1), F)$. Note the similarity with $\Der_F(X,M) = \Map (\RuSpec(A \oP M),F)$. Now by Proposition 1.4.1.6 of \cite{TV4}, for $F \in \St(\cC,\tau)$, $A \in \cA$, $x: X= \bR \uSpecA \rarr F$, supposing $F$ has a cotangent complex $\LFx$ at $x$, it follows we have natural isomorphisms:
\beq
\RuHom_{\AffCtildetau}(X,TF) \cong \Map_{\Sp(\AMod)}(\LFx,A) \cong \Map_{\Sp(\AMod)}(A,\TFx) \nonumber
\eeq
thereby showing that $\TFx$ is the analog among spectra of the tangent stack, which justifies calling $\LFx$ the cotangent complex of $F$ at $x$. To see that we have those isomorphisms, it suffices to write:
\begin{align}
	\RuHom(X,TF) &= \RuHom(\bR \uSpecA, \RuHom(\bD_{\eps},F)) \nonumber \\
	&\cong \RuHom( \bR \uSpecA \times^h \bD_{\eps},F) \nonumber \\
	&\cong \RuHom( \bR \uSpec(A[\eps]),F) \nonumber \\
	&= \RuHom(\RuSpec(A \oP A),F) \nonumber \\
	&=\RuHom(X[A],F) \nonumber \\
	&=\Der_F(X,A) \nonumber \\
	&\cong \bR \uh_s^{\LFx}(A) = \Map(\LFx,S_A(A))\nonumber 
\end{align}
Since $\skMod$ has a fully faithful suspension functor, by lemma 1.2.11.2 of \cite{TV4} the stabilization functor $S_A$ is fully faithful as well, so we just write $A$ for $S_A(A)$ in what follows: 
\begin{align}
	\RuHom(X,TF) &\cong \Map(\LFx,A) \nonumber \\
	&\cong \Map(\LFx \oTAL A,A) \nonumber \\
	&\cong \Map(A, \RuHom(\LFx,A)) \nonumber \\
	&=\Map(A,\TFx) \nonumber
\end{align}

\subsubsection{Convergence statement}
The last claim of that representability theorem states that for any $A \in \skAlg$, $\bR F(A) \xrarr{\cong} \holim_k \bR F(A_{\leq k})$. As pointed out in \cite{TV4}, any $A \in \skAlg$ has a functorial Postnikov tower of the form:
\beq
A \rarr \cdots \rarr A_{\leq n} \rarr A_{\leq n-1} \rarr \cdots \rarr A_{\leq 0} = \pi_0(A) \nonumber
\eeq
where $A_{\leq k}$ is $k$-truncated, the morphisms $A \rarr A_{\leq k}$ induce isomorphisms on fundamental groups in degree $l \leq k$, and most importantly for us, this Postnikov tower is unique up to equivalence. In other terms one can define the \textbf{$k$-truncation} of $A$ to be given by $A_{\leq k}$ in its Postnikov tower. In this manner one can write:
\beq
F(A_{\leq k}) = F(tr_k(A)) = F \circ tr_k(A) := F_k(A) \nonumber
\eeq
where $tr_k$ is the degree $k$-truncation functor, which exists by virtue of the above discussion, and $F_k$ is defined as above. We call it the degree $k$-approximation to the functor $F$.

\section{The bridge} \label{The bridge}
This section aims to show that $\rho$-analyticity for Goodwillie translates into a stack being $n$-geometric, plus some added conditions.
\subsection{Similitudes in towers}
$F \in \kDAfftildeet$ maps the Postnikov tower of $A \in \skAlg$ to:
\beq
F(A) \rarr \cdots F(A_{\leq n}) \rarr F(A_{\leq n-1}) \rarr \cdots \rarr F(A_0) = F(\pi_0(A)) \nonumber
\eeq
which one can rewrite like so:
\beq
F(A) \rarr \cdots F_n(A) \rarr F_{n-1}(A) \rarr \cdots \rarr F_0(A) = F(\pi_0(A)) \nonumber
\eeq
In $\Ho(\SetD)$, this reads:
\beq
\bR F(A) \rarr \cdots \bR F_n(A) \rarr \bR F_{n-1}(A) \rarr \cdots \rarr \bR F_0(A) = F(\pi_0(A)) \nonumber
\eeq
since $F(A) \cong \bR F(A)$ in $\Ho(\SetD)$, where $\bR F(A) = \bR_{\tau} \uHom(\bR \uSpecA, F)$, $\bR_{\tau} \uHom$ the derived simplicial Hom for $\kDAfftildeet$, which, because $F$ is a stack, can also be written using $\RuHom$, the derived simplicial hom of $\kDAffhat$.\\

\newpage

Now if $F$ is a $n$-geometric $D^{-}$-stack, we have $\bR F(A) \cong \holim_n \bR F_n(A)$, which one can view, in the spirit of functor Calculus, as a convergence of the above tower for $F$ at the point $A$. \\

Another hint that there may be a link between those Taylor towers, whether they be for a geometric stack, or for a homotopy functor, comes from the following considerations. As pointed out in \cite{TV4}, the homotopy fiber of the morphism $A_{\leq n} \rarr A_{\leq n-1}$ in the Postnikov tower of $A \in \skAlg$ is isomorphic to $S^n \oT_k \pi_n A$ in $\Ho(\skMod)$, where the concept of homotopy fiber is as defined in \cite{Hi}; for a functorial factorization $A_{\leq n} \rarr E(A_{\leq n}) \rarr A_{n-1}$, a trivial cofibration followed by a fibration, the homotopy fiber of $A_{\leq n} \rarr A_{\leq n-1}$ is defined by $P_n = * \times_{A_{n-1}} E(A_{\leq n})$, thus the following diagram in particular is a homotopy cartesian square:
\beq
\xymatrix{
P_n \ar[d] \ar[r] & \Aln \ar[d] \\
\text{*} \ar[r] & \Alnmone
} \nonumber
\eeq
Suppose now $F \in \kDAfftildeet$ preserves homotopy cartesian squares. It follows:

\beq
\xymatrix{
	 F(P_n) \ar[d] \ar[r] &  F(A_{\leq n}) \ar[d] \\
	 F(*) \ar[r] &  F(A_{\leq n-1})
} \nonumber
\eeq
is a homotopy cartesian square as well, that is:
\beq
F(P_n) \xrarr{\simeq} F(*) \times_{F(\Alnmone)} E(F \Aln) \nonumber
\eeq
Let $\gamma: \SetD \rarr \Ho(\SetD)$ be the canonical projection. It follows:
\beq
\gamma F(P_n) \xrarr{\cong} \gamma F(*) \times_{\gamma F(\Alnmone)} \gamma E(F \Aln) \nonumber
\eeq
Suppose $F$ is reduced: $F(*) = *$, so that $\gamma F(*) = \gamma * = *$. This reads:
\beq
\gamma F(P_n) \xrarr{\cong} * \times_{\gamma F(\Alnmone)} \gamma E(F \Aln) \nonumber
\eeq
Now $\gamma F: \skAlg \rarr \Ho(\SetD)$ factors uniquely through $\Ho(\skAlg)$: $\Ho(\gamma F) : \Ho (\skAlg) \rarr \Ho(\SetD)$, which we denote again by $\gamma F$. Using the fact that the homotopy fiber $P_n$ of $\Aln \rarr \Alnmone$ is isomorphic to $S^n \oT_k \pi_n A$, which we denote by $\Omega M$, it follows:
\beq
\gamma F(\Omega M) \cong * \times_{\gamma F(\Alnmone)} \gamma E( F \Aln) \nonumber
\eeq
We have an equivalence $F(A) \xrarr{\simeq} \bR F(A)$ in $\SetD$ for $A \in \skAlg$ since $F$ is a stack (hence an isomorphism $\gamma F(A) \cong \bR F(A)$ in $\Ho(\SetD)$), so that we have an equivalence $E (FA) \simeq E(\bR F A)$ by the 2-3 property, and consequently an isomorphism $\gamma E (FA) \cong \gamma E( \bR FA)$ in $\Ho(\SetD)$. It follows:  
\beq
\bR F(\Omega M) \cong * \times_{\bR F(\Alnmone)} \gamma E(\bR F(\Aln)) \nonumber
\eeq
Since we work in $\Ho(\SetD)$, we can drop $\gamma$. This now reads:
\beq
\bR F (\Omega M) \cong * \times_{\bR F(\Alnmone)} E(\bR F(\Aln)) \nonumber
\eeq
Suppose now $F$ is a monoidal functor: $F(A) \oT_{\SetD} F(B) \xrarr{\cong} F(A \oT_k B)$. Using the fact that $\Omega M = S^n \oT_k \pi_n A$, this reads:
\beq
\bR F(S^n) \oT_{\SetD} \bR F(\pi_n A) \cong * \times_{\bR F(\Alnmone)} E (\bR F( \Aln)) \nonumber
\eeq
which means $\bR F(S^n) \oT \bR F(\pi_n A)$ is the homotopy fiber of $\bR F_n(A) \rarr \bR F_{n-1}(A)$ in $\Ho(\SetD)$.\\

Now recall from \cite{Go3} that for $F$ a homotopy functor, the homotopy fiber of $P_n F(X) \rarr P_{n-1}F(X)$, denoted $D_n F(X)$, is, for spectra valued functors, of the form $(C_n \wedge (X^{\wedge n}))_{h \Sigma_n}$ for some spectrum $C_n$ with a $\Sigma_n$ action. This is surprisingly similar to what we have for our stack. The correspondence is given by associating $C_n$, which Goodwillie calls the $n$-th derivative of $F$ at $x$ and which he denotes by $\del^{(n)}F(x)$, with $F(S^n)$, and $X^{\wedge n}$ with $F(\pi_n A)$. The corresponding action of the symmetric group $\Sigma_n$ is given by permuting the entries of $[n]$ in $\Dn = \Hom(-,[n])$. Note that here $X$ is regarded as a space over the one point space $*$, so likewise $A \in \skAlg$ is regarded as an algebra over 1. Thus in addition to having convergence of Taylor towers for $\rho$-analytic functors and for $n$-geometric $D^{-}$-stacks, incremental morphisms in such towers have homotopy fibers that are formally the same if in addition we focus on those stacks that preserve homotopy cartesian diagrams, are reduced and monoidal.\\

\newpage

\subsection{Analyticity and obstruction theory: a study}
Now that we have seen we have comparable phenomena at the level of towers, we study this correspondence more carefully. We will now see that what makes $F$ a $n$-geometric stack in $\DStk$ is very close indeed to the conditions that Goodwillie put on a homotopy functor to have convergence of its Taylor tower. Coming back to the condition that $F$ preserves homotopy cartesian squares, observe that this is a scaled down condition on a stack $F$ relative to what Goodwillie asked his functors to satisfy for convergence, namely that they map strongly co-cartesian squares to cartesian squares. It turns out this is the only such condition of the sort we will impose on our stacks, and for this reason we refer to it as the \textbf{op-excisive} condition, for obvious reasons. \\

In the representability theorem for stacks Theorem C.0.9 we are considering, $F$ has an obstruction theory. We will show this is morally the same for stacks as asking a homotopy functor to be analytic. For a stack $F$ to have an obstruction theory necessitates that it be infinitesimally cartesian and that it has a global cotangent complex. Recall what infinitesimally cartesian means: given the homotopy cartesian square
\beq
\xymatrix{
	A \oP_d \Omega M \ar[d] \ar[r] & A \ar[d] \\
	A \ar[r] &A \oP M
} \nonumber
\eeq
the following square is homotopy cartesian as well:
\beq
\xymatrix{
	\bR F(A \oP_d \Omega M) \ar[d] \ar[r] & \bR F(A) \ar[d] \\
	\bR F(A) \ar[r] &\bR F(A \oP M)
} \nonumber
\eeq
This is automatic if $F$ is op-excisive. Thus asking that an op-excisive stack $F$ has an obstruction theory boils down to showing it has a cotangent complex, which we will show corresponds in a sense to the stably excisive condition, needed for the convergence of Taylor towers of homotopy functors.\\

To ask that $F$ has a cotangent complex means we have a $(-n)$-connective stable module $\LFx$ such that $\Der_F(X,-) \cong \bR \uh_s^{\LFx}$ as we have seen, where $X = \bR \uSpecA$ and $x: X \rarr F$. The derivation from $F$ functor $\Der_F(X,-)$ has a lift to a functor $\AMod \rarr \SetD$, which to $M$ associates the homotopy fiber at $x \in F(A)$ of the map $F(A \oP M) \rarr F(A)$ for any $M \in \AMod$. This is the image under $F$ of the map $A \oP M \rarr A$. Collectively those maps can be represented by a functor $\cX : \AMod \rarr \skAlg$, which to $M$ associates $\cX(M) = A \oP M$, in such a manner that $A \oP M \rarr A$ reads $\cX(M) \rarr \cX(\eset)$. Now observe that for Goodwillie, being $\rho$-analytic for a homotopy functor $F$ means that it is $n$-stably excisive, meaning it satisfies $E(n,c,K)$, which we reproduce here: if $\cX : \cP(S) \rarr \cC$ is a strongly co-cartesian $(n+1)$-cube such that $\cX(\eset)\rarr \cX(s)$ is $k_s$-connected for all $s \in S$, with $k_s \geq K$, then $F(\cX)$ is $(-c + \sum k_s)$-cartesian. Thus on the one hand one considers maps of the form $\cX(\eset) \rarr \cX(s)$ in the calculus of functors, while for stacks one considers maps of the form $A \oP M \rarr A \oP \eset = A$ instead, which read $\cX(M) \rarr \cX(\eset)$. The fact that we use ``wrong" direction maps is explained by the fact that we ask that $F$ preserves homotopy cartesian diagrams, as opposed to mapping homotopy co-cartesian diagrams to homotopy cartesian diagrams in the case of functor calculus. However, the notion of cubical diagram here is moot for stacks; $\AMod$, the indexing set for $\cX$, is not finite. No need to work with strongly co-cartesian diagrams either, since our stack $F$ is assumed to be op-excisive. In the last section of this work we will see that picking $\cX = A \oP -$ is not simply an ad hoc identification; the convergence of the Taylor tower of a homotopy functor $F$ at some space $X$ over $Y$ hinges on the choice of $\cX = X *_Y -$ as cubical diagram, which we identify with $A \oP -$ in Algebraic Geometry.\\

 Note that whether it be a Taylor tower, or a cotangent complex, what we have is very much deformation theoretic, thus we will put an emphasis on connectivity rather than on connectedness. This choice is further reinforced by the fact that if $\cX(\eset) \rarr \cX(s)$ is $k_s$-connected, $F(\cX)$ is $(-c + \sum k_s)$-cartesian, i.e. by Def. 1.3 of \cite{Go2}, the homotopy fibers of $F \circ \cX(\eset) \simeq \holim F \circ \cX \rarr \holim( F \circ X|_{\cP_0(S)})$ are $(-c + \sum k_s -1)$-connected. Collectively we have preserved the connectedness of each map. We will show that the same is true of connectivity in the case of stacks.\\

Now that we have defined diagrams for stacks, we study more carefully what the assumptions and claims made on diagrams in the calculus of functors translate to in the language of derived stacks. For Goodwillie, $\cX(\eset) \rarr \cX(s)$ being $k_s$-connected means its homotopy fibers are $(k_s-1)$-connected, and here we are considering homotopy (limit) fibers. In the case of stacks, $\cX = A \oP -$, the homotopy fibers of morphisms $A \oP M \rarr A$ depend on $M$. Thus one defines the connectivity of the map $A \oP M \rarr A$ to be the one for $M$, where we say $M$ is $(-p)$-connective if it is of the form $M = \Omega^p M'$ for some $M' \in \AMod$ $-(p-1)$-connective. Because we have such maps, they functorially map to $F(A \oP M) \rarr F(A)$, the collection of which forms $F(\cX)$. For Goodwillie, $F \cX$ being $(-c + \sum k_s)$-cartesian means the homotopy fibers of $F \cX (\eset) \rarr \holim(F|_{\cP_0(S)})$ are $(-c + \sum k_s -1)$-connected. Instead of looking at the fibers of the map $\holim F(A \oP M)  \rarr F(A)$ to have a statement mirroring the one made by Goodwillie for $F(\cX)$, we consider the homotopy fibers of each individual map $F(A \oP M) \rarr F(A)$ instead, all of which, as pointed out above, are functorially represented by $\Der_F(X,-)$, with $X = \bR \uSpecA$. Thus for the calculus of functors, the $n$-stably excisive condition for functors $F$ focuses on $F \cX$, while for stacks one looks at $\Der_F(X,-)$ instead. For stacks, the observation made above regarding the preservation of the connectedness of objects in the calculus of functors(captured by the presence of $\sum k_s$) should translate into the preservation of connectivity as we argued. Indeed, for $F$ a $n$-geometric $D^{-}$-stack, $\Der_F(X,\Omega M) \cong \Omega \Der_F(X,M)$, thus if $A \oP M \rarr A$ has connectivity $p$, meaning $M$ is a $p$-loop object, so does $\Der_F(X,M)$, whose lift is the homotopy fiber of $F(A \oP M) \rarr F(A)$, thus $F$ preserves the connectivity of each map. To address the occurrence of $-c$ in the $n$-stably excisive condition, observe that $\Der_F(X,-)$ is co-represented by $\LFx$, which is connective. \\

Regarding the connectivity of $\LFx$, by Lemma 1.4.3.10 of \cite{TV4} a $n$-geometric $D^{-}$-stack $F$ is infinitesimally cartesian, so by the proof of Proposition 1.4.2.7 of the same reference, it satisfies the conditions of Proposition 1.4.1.11, from which it follows that $F$ has a global cotangent complex, which is $(-n)$-connective, with the same index $n$. \\ 

Thus we have shown that for $F$ an op-excisive (reduced, monoidal) $n$-geometric $D^{-}$-stack, $\rho$-analyticity translates into $F$ having an obstruction theory.\\

\section{Parallel results} \label{Parallel results}
\subsection{Polynomial approximations}
At this point we can draw various comparisons between both theories. On the one hand for $F: \cC_Y \rarr \cD$ a homotopy functor, $X \in \cC_Y$, we have:
\begin{align}
	P_nF(X) &= \hocolim_k T_n^k F(X) \nonumber \\
	&\cong \hocolim_k( \holim_{(U_1, \cdots, U_k) \in \cP_0(\unpone)^k} F(X *_Y(U_1 * \cdots * U_k))) \nonumber
\end{align}
with $X *_Y U = \hocolim(X \larr X \times U \rarr Y \times U)$. We compare this with $\bR F_n(A) = \bR F(\Aln)$, $\pi_i (\Aln) = 0$ if $i > n$. From \cite{GoJa}, $\Aln = \DdA_{/\sim_n}$ where $\alpha \sim_n \beta$ for $\alpha, \beta: \Delta^q \rarr A$ if their restrictions to $\skn(\Delta^q)$ coincide. One can write:
\beq
	\Aln = \cup_{q \geq 0} \coprod_{\alpha \in \DqdA} [\alpha] \nonumber 
\eeq
with $[\alpha] = \{ \beta \in \DqdA \, | \, \alpha \sim_n \beta \}$. If $\alpha,\beta \in \DqdA$, define $\alpha \star_A^q \beta$ by $\alpha \cup \beta$ if $\alpha \sim_n \beta$, $\eset$ otherwise. With this if $F \in \kDAfftildeet$:
\begin{align}
	\bR F_n(A) &= \bR F(\cup_{q \geq 0} \coprod_{\alpha} \cup_{k \geq 0} \alpha \star_A^q (\beta_1 \star_A^q \cdots \star_A^q \beta_k )) \nonumber \\
	&=\bR F(\cup_{k \geq 0} \cup_{q \geq 0} \coprod_{\alpha \in \DqdA} \bigcup_{\substack{\beta_l \sim_n \alpha \\ 1 \leq l \leq k}}\alpha \star_A^q ( \beta_1 \star_A^q \cdots \star_A^q \beta_k )) \nonumber \\
	&=\bR F(\cup_{k \geq 0} \cup_{q \geq 0} \bigcup_{\substack{\beta_l \in \DqdA \\ 1 \leq l \leq k}} (\DqdA) \star_A^q (\beta_1 \star_A^q \cdots \star_A^q \beta_k)) \nonumber \\
	&=\bR F( \colim_k \cup_{q \geq 0} \colim_{(\beta_1, \cdots, \beta_k) \in \text{Sk}_{0,q}(n)^k} \DqdA \star_A^q (\beta_1 \star_A^q \cdots \star_A^q \beta_k)) \nonumber 
\end{align}
where $\text{Sk}_{0,q}(n) = \{ \beta \in (\DqdA) - \eset \, | \, \alpha \sim_n \beta \}$. With $\DqdA \cong A_q$, this reads:
\beq
\bR F_n(A) = \bR F( \colim_k \cup_{q \geq 0} \colim_{(\beta_1, \cdots, \beta_k) \in \text{Sk}_{0,q}(n)^k}A_q \star_A^q (\beta_1 \star_A^q \cdots \star_A^q \beta_k)) \nonumber
\eeq
Denoting $\text{Sk}_0(n)^k = \cup_{q \geq 0} \text{Sk}_{0,q}(n)^k$, write:
\beq
\colim_{(\beta_1, \cdots, \beta_k) \in \text{Sk}_0(n)^k} A \star_A (\beta_1 \star_A \cdots \star_A \beta_k) = \cup_{q \geq 0} \colim_{(\beta_1, \cdots, \beta_k) \in \text{Sk}_{0,q}(n)^k} A_q \star_A^q (\beta_1 \star_A^q \cdots \star_A^q \beta_k) \nonumber
\eeq
where $\star_A$ is a levelwise product, so that:
\beq
\bR F_n(A) = \bR F( \colim_k \colim_{(\beta_1, \cdots , \beta_k) \in \text{Sk}_0(n)^k} A \star_A (\beta_1 \star_A \cdots \star_A \beta_k)) \nonumber
\eeq
which one can compare with:
\beq
P_nF(X) \cong \hocolim_k \holim_{(U_1, \cdots, U_k) \in \cP_0(\unpone)^k} F(X *_Y(U_1 * \cdots * U_k)) \nonumber
\eeq

\subsection{Pointwise approximations}
We now consider Proposition 1.6 of \cite{Go3}, which states in particular that if $F$ and $G$ are two stably $n$-excisive homotopy functors, $u:F \Rarr G$, then $F$ and $G$ satisfy the following condition: there are constants $c$ and $K$ such that $\forall \, k \geq K$, $\forall \, X \in \cC_Y$ for which $X \rarr Y$ is $k$-connected, the map $F(X) \rarr G(X)$ is $(-c + (n+1)k)$-connected. The fact that $F(X) \rarr G(X)$ is $(-c + (n+1)k)$-connected $\forall X$ shows that there is an intrinsic connective object to $u$. Now being $n$-stably excisive is a property of $\rho$-analytic functors, which we associate with $n$-geometric $D^{-}$-stacks in the theory of stacks. Thus to show for stacks a result similar to Proposition 1.6, we start with two $n$-geometric $D^{-}$-stacks $F$ and $G$ and a morphism of stacks $f: F \rarr G$. We consider maps of the form $f_A: F(A) \rarr G(A)$ for $A \in \skAlg$, and we aim to show all such maps are inherently connective with the same connectivity in a sense to be precised. Since we work deformation theoretically, we will lift those maps to all diagrams of the form:
\beq
\xymatrix{
	F(A \oP M) \ar[d] \ar@{.>}[r] & G(A \oP M) \ar[d] \\
	F(A) \ar[r] & G(A)
}\label{GCGSCD}
\eeq
indexed by $M \in \AMod$, whose individual homotopy fibers are functorially represented by derivation functors. In lieu of considering \eqref{GCGSCD}, we therefore focus on the following natural transformation:
\beq
\Der_F(X,-) \Rarr \Der_G(X,-) \nonumber
\eeq
with $X = \bR \uSpec A$, in the same manner that we had $u_X: F(X) \rarr G(X)$ in Proposition 1.6 above. To be precise, for $x: X=\bR \uSpecA \rarr F$ an $A$-point of $F$, we have an $A$-point $ f \circ x: \bR \uSpecA \rarr G$ of $G$, and the natural transformation above appears as an induced map:
\beq
f_* = df: \Der_F(X,-) \Rarr \Der_G(X,-) \nonumber
\eeq
where $\Der_G(X,-)$ is taken at the point $f(x)$. Since in Proposition 1.6 $u_X: F(X) \rarr G(X)$ is $(-c + (n+1)k)$-connected, its homotopy fibers are $(-c + (n+1)k-1)$-connected. This motivates the introduction of $\Der_{F/G}(X,-)$, the homotopy fiber of $df$ (\cite{TV4}), defined by:
\beq
\Der_{F/G}(X,M) = \Map_{X/\kDAfftildeet/G}(X[M],F) \nonumber
\eeq
and we want to show it has an intrinsic connective object. To that end, observe that $f$ has a cotangent complex, which is connective, as we now explain. By Def. 1.4.1.14 of \cite{TV4}, $f: F \rarr G$ has a \textbf{relative cotangent complex at $x$} relative to $\skAlg$ if there is some $(-n)$-connective stable $A$-module $\LFGx \in \Ho(\Sp(\AMod))$ and an isomorphism $\Der_{F/G}(X,-) \cong \bR \uh_s^{\LFGx}$. If this is true, the homotopy fiber of $\Der_F(X,-) \Rarr \Der_G(X,-)$ would be connective. However note that this is only for $A$ fixed, whereas we have $F(X) \rarr G(X)$ $(-c + (n+1)k)$-connected for all $X \in \cC_Y$ in the calculus of functors. This means the relative cotangent complexes we just introduced must be found for all points $A \in \skAlg$ and exhibit the same connectivity. Consequently we need Def. 1.4.1.15 of \cite{TV4}: we say $f$ has a \textbf{cotangent complex relative to $(\skMod,\skMod,\skAlg)$} if $\forall \, A \in \skAlg$, $\forall \, x: \bR \uSpecA \rarr F$, $f$ has a cotangent complex $\LFGx$ at $x$, and additionally $\forall u:A \rarr B$ in $\skAlg$, any morphism:
\beq
\xymatrix{
	\bR \uSpec B \ar[dr]_y \ar[rr] &&\bR \uSpecA \ar[dl]^x \\
	&F
} \nonumber
\eeq
in $\Ho(\kDAfftildeet)$, induces an isomorphism $\LFGx \oTAL B \xrarr{\cong} \LFGy$ in $\Ho(\Sp(\BMod))$, a sort of coherence condition. But by lemma 1.4.1.16 of \cite{TV4}, if both $F$ and $G$ have a cotangent complex, so does $f$. However by Theorem C.0.9, if both $F$ and $G$ are $n$-geometric $D^{-}$-stacks, they have an obstruction theory, in particular they have a cotangent complex, hence so does $f$. In other terms, lifting each morphism $F(A) \rarr G(A)$ to the level of derivations, we obtain a natural transformation $\Der_F(X,-) \Rarr \Der_G(X,-)$ with $X = \bR \uSpecA$, whose homotopy fiber $\Der_{F/G}(X,-)$ at $x$ has a (connective) cotangent complex, and $\forall u: A \rarr B$ in $\skAlg$, $x: \bR \uSpec A \rarr F$, $y: \bR \uSpec B \rarr F$, those cotangent complexes are compatible in the following sense: $\LFGx \oTAL B \xrarr{\cong} \bL_{F/G,y}$, which is our analog for stacks of Proposition 1.6 for homotopy functors.

\subsection{Reconstruction theorem}
We have the following result from \cite{AC}, Theorem 2.9, a reconstruction theorem for the Taylor tower of a homotopy functor $F: \Sp \rarr \Sp$, which states that there is a natural homotopy cartesian square of the form:
\beq
\xymatrix{
	P_nF(X) \ar[d] \ar[r] & (\del_n F \wedge X^{\wedge n})^{h \Sigma_n} \ar[d] \\
	P_{n-1}F(X) \ar[r] & (\del_n F \wedge X^{\wedge n})^{t \Sigma_n}
} \nonumber
\eeq
where if $X$ is a spectrum with an action of a finite group $G$, $X^G$ the spectrum of $G$-invariants of $X$, $X^{tG}$ is the cofiber of the norm map $X/G \rarr X^G$. For our part, by Lemma 2.2.1.1. of \cite{TV4}, there is a unique derivation $d_n \in \pi_0( \Der_k(\Alnmone, S^{n+1} \oT_k \pi_n A))$ for which we have an isomorphism in $\Ho(\skAlg_{/\Alnmone})$ between the natural projection $\Alnmone \oP_{d_n} (S^n \oT_k \pi_n A) \rarr \Alnmone$ and $\Aln \rarr \Alnmone$. For this particular derivation, we have a homotopy cartesian square:
\beq
\xymatrix{
	\Alnmone \oP_{d_n} \Omega(S^{n+1} \oT_k \pi_n A) \ar[d] \ar[r] & \Alnmone \ar[d]^{d_n} \\
	\Alnmone \ar[r]_-s & \Alnmone \oP (S^{n+1} \oT_k \pi_n A)
} \nonumber
\eeq
Let $M = S^{n+1} \oT_k \pi_n A$. The reasoning now is the same as for showing homotopy fibers of connecting maps in towers are similar in form. Because $F$ is op-excisive, we have a homotopy cartesian square:
\beq
\xymatrix{
	F(\Alnmone \oP_{d_n} \Omega M) \ar[d] \ar[r] & F(\Alnmone) \ar[d] \\
	F(\Alnmone) \ar[r] & F(\Alnmone \oP M)
} \nonumber
\eeq
Thus $F(\Alnmone \oP_{d_n} \Omega M) \xrarr{\simeq} F \Alnmone \times_{ F(\Alnmone \oP M)} E (F \Alnmone)$. By Proposition 1.3.3.8 of \cite{Hi}, for $A,B, C$ objects of some right proper model category $\cC$, $A \times_C E(B)$ is weakly equivalent to $E(A) \times_C E(B) = A \times^h_C B$, hence an isomorphism in $\Ho(\cC)$. It follows that $\gamma F(\Alnmone \oP_{d_n} \Omega M) \xrarr{\cong} \gamma F \Alnmone \times^h_{ \gamma F(\Alnmone \oP M)} \gamma F \Alnmone$. Now $\gamma F: \skAlg \rarr \Ho(\SetD)$ has a unique factoring map $\Ho(\gamma F): \Ho(\skAlg) \rarr \Ho(\SetD)$ also denoted $\gamma F$ for simplicity. Further, using the fact that $\Omega M \cong S^n \oT_k \pi_n A$, and the fact that the map $\Alnmone \oP_{d_n} S^n \oT_k \pi_n A \rarr \Alnmone$ is isomorphic to $\Aln \rarr \Alnmone$, it follows that we have the following isomorphism in $\Ho(\SetD)$:
\beq
\gamma F(\Aln) \xrarr{\cong} \gamma F(\Alnmone) \times^h_{\gamma F(\Alnmone \oP M)} \gamma F \Alnmone \nonumber
\eeq
Finally for $A \in \skAlg$, $F(A) \xrarr{\simeq} \bR F(A)$, hence an isomorphism $\gamma F(A) \xrarr{\cong} \bR F(A)$, from which it follows that we have a homotopy pullback square in $\Ho(\SetD)$:
\beq
\xymatrix{
	\bR F(\Aln) \ar[d] \ar[r] & \bR F(\Alnmone) \ar[d] \\
	\bR F(\Alnmone) \ar[r] & \bR F(\Alnmone \oP (S^{n+1} \oT_k \pi_n A))
}\nonumber
\eeq
which provides a reconstruction as well of the tower for $\bR F(A)$.

\section{Taking stock - General theory} \label{Taking stock}

\subsection{Convergence of the Taylor tower}
Going over the proof of Theorem 1.13 of \cite{Go3}, let $F$ be $\rho$-analytic, so that it is $n$-stably excisive. In the definition of $E(n,c,K)$,  we pick $\cX = X *_Y -$ as a cubical diagram  from $\cP(\unpone)$ into $\cC$. If $s \in \unpone$, $\cX(\eset) \rarr \cX(s)$ reads $X \rarr X *_Y s$. Now for the convergence of the Taylor tower, we suppose $X \rarr Y$ is $k$-connected. By Remark 1.1 of \cite{Go3} it follows that $X *_Y s \rarr Y$ is $(k+1)$-connected, and by by Proposition 1.5 of \cite{Go2} that $X \rarr X *_Y s$ is $k$-connected $\forall s \in \unpone$, so that $F \cX$ is $(-c + (n+1)k)$-cartesian since $F$ satisfies $E(n,c,K)$, which exactly means that $F\cX(\eset) \rarr \holim(F\cX|_{\cP_0(\unpone)})$ is $(-c + (n+1)k)$-connected. This map reads $F(X) \rarr \holim_{U \in \cP_0(\unpone)}F(X *_Y U) = T_nF(X)$. By Def. 1.2, this is referred to as saying that $F \rarr T_n F$ satisfies $\cO(n,c,K)$. In the statement of Theorem 1.13, $X \rarr Y$ is $(\rho + 1)$-connected, where $k = \rho + 1$ is that particular value that allows us to have convergence of the Taylor tower. Thus a key step in proving the convergence of the Taylor tower for $F(X)$, with $F$ $\rho$-analytic, $X \rarr Y$ $(\rho+1)$-connected, consists in first picking $\cX = X *_Y -$ from which we have $F \rarr T_nF$ satisfying $\cO(n,c,K)$ for some $c$ and $K$, which is interpreted as saying that $F$ and $T_nF$ agree to order $n$. Aside from the fact that $T_nF$ will be useful in constructing $P_nF$, this degree $n$ approximation is essential to prove the convergence of the Taylor tower for $F(X)$.\\

Note that it is at this point that we took the analogous diagram for stacks to be represented by $A \oP M \rarr A$, which amounts to identifying $X *_Y -$ with $A \oP -$.
\subsection{Convergence of the tower for stacks}
Presently, we focus on the convergence of the tower involving $\bR F(\Aln)$. Essentially we have $\bR F(A) \cong \holim_n \bR F(\Aln)$ by virtue of the fact that $\Aln$ can be expressed as a square zero extension. More specifically, following the proof of Theorem C.0.9 of \cite{TV4} to show that we have such an isomorphism in the homotopy category of simplicial sets, we show we have an equivalence in $\SetD$, or equivalently that for some $x \in \pi_0(\holim_n \bR F(\Aln))$, the homotopy fiber of the natural map $\bR F(A) \rarr \holim_n \bR F(\Aln)$ is contractible. To achieve this one uses induction. If $x_n \in \bR F(\Aln)$ is the projection of $x$, one argues there is a representable $D^{-}$-stack $U$ along with a smooth morphism $U \rarr F$ such that one has $y_n \in \pi_0(\bR U(\Aln))$ mapping to $x_n$. Because $\Aln$ is a square zero extension of $\Alnmone$, one finds that the morphism $\bR U(A_{\leq n+1}) \rarr \bR U(\Aln) \times^h_{\bR F(\Aln)} \bR F(A_{\leq n+1})$ is surjective, which allows us to lift $y_n$ to some element $y_{n+1} \in \pi_0(\bR U(A_{\leq n+1})$, and from there one obtains an element $y \in \pi_0(\holim \bR U(\Aln))$ mapping to $x$. Thus we get a commutative diagram:
\beq
\xymatrix{
	\bR U(A) \ar[d] \ar[r] & \holim \bR U(\Aln) \ar[d] \\
	\bR F(A) \ar[r] & \holim \bR F(\Aln) 
}\nonumber
\eeq
By induction this is a homotopy cartesian square, and the top map is an equivalence, which allows us to conclude. All the results necessary to conduct this proof are based on the fact that $F$ has an obstruction theory, and the fact that the Postnikov tower for $A \in \skAlg$ is constructed from square zero extensions, thus the underlying reasoning is purely deformation theoretic.\\

To emphasize this last point, from Lemma 2.2.1.1, there is a unique derivation $d_n \in \pi_0 (\Der_k(\Alnmone, S^{n+1} \oT_k \pi_n A))$ such that $\Aln \rarr \Alnmone$ is isomorphic to $\Alnmone \oP_{d_n} S^n \oT_k \pi_n A \rarr \Alnmone$. Thus the tower for $F_n(A) = F(\Aln)$ reads:
\beq
\cdots \rarr \bR F(\Alnmone \oP_{d_n} (S^n \oT_k \pi_n A)) \rarr \bR F(\Alnmone) \rarr \cdots \nonumber
\eeq
It follows that starting from an element of $\pi_0 F(A)$, to climb up in the tower for $\bR F(A)$ one has to iterate moves of the form $\bR F(A) \rarr \bR F(A \oP_d \Omega M)$ for $M \in \AMod_1$. Using the fact that $\bR F(A) \cong \RuHom( \bR \uSpec A, F)$, this amounts to having a lift of $x$ to $x'$ in the following diagram:
\beq
\xymatrix{
	\bR \uSpec A \ar[dr]_x \ar[rr] && \bR \uSpec (A \oP_d \Omega M) \ar[dl]^{x'} \\
	&F
}\nonumber
\eeq
However by Proposition 1.4.2.5, we have such a lift if we have a vanishing obstruction class in $\Map_{\Sp(\AMod)}(\LFx,M)$. Thus we obtain a construction that is inherently deformation theoretic.\\

\subsection{Abstraction}
For Goodwillie calculus, we used the cubical diagram $X *_Y -$ to obtain that $F(X) \rarr \holim F(X *_Y U)$ is $(-c + (n+1)k)$-connected. This is a pivotal statement for proving convergence of the Taylor tower of $F(X)$. One can note that assuming $F$ is $\rho$-analytic, if $X$ is close enough to $Y$ in the sense that $X \rarr Y$ is $\rho +1$-connected, then we have convergence of the Taylor tower if $n$ goes to infinity in $F(X) \rarr P_nF(X)$. Essentially, we have convergence in a neighborhood of $Y$.\\

\newpage

For stacks we took our analogous diagram to be given by $\cX = A \oP -$, and in lieu of dealing with $F\cX$, we considered all maps $F(A \oP M) \rarr F(A)$ and argued that functorially their homotopy fibers are a lift of the functor $\Der_F(X,-)$, which because $F$ has an obstruction theory, can be co-represented by $\LFx$. Having a cotangent complex $\LFx$ allows us to use obstruction theory arguments to eventually prove the convergence of $\bR F(\Aln)$ to $\bR F(A)$. \\

So far both theories are built upon the choice of a deformation theoretic diagram, whether it be $X *_Y -$ or $A \oP -$, and a pivotal result, $F \rarr T_nF$ agreeing to order $n$ in Goodwillie's calculus, and the existence of $\LFx$ for stacks. We also had a notion of neighborhood for the calculus of functors. We now investigate such a notion for stacks.\\

In order to do so, observe that convergence of the Postnikov tower for $\bR F(A)$ hinges on Lemma C.0.10, in which proof by using a local argument, one argues one can find a representable $D^{-}$-stack $U$ along with a smooth morphism $U \rarr F$, in such a manner that we have the following summarizing diagram:
\beq
\xymatrix{
	&\bR U(A) \ar[dl]_{\cong} \ar[r]^{fctrial} & \bR F(A) \ar[dr]^{nat.map} \\
	\holim \bR U(\Aln) \ar@{.>}[dr] \ar@{.>}[ddr] \ar@{.>}[rrr] &&& \holim \bR F(\Aln) \ar[dl] \ar[ddl] \\
	&\bR U(A_{\leq n+1}) \ar@{.>}[d] \ar[r] & \bR F(A_{\leq n+1}) \ar[d] \\
	&\bR U(\Aln) \ar[r] & \bR F(\Aln) 
} \nonumber
\eeq
The argument goes as follows: one can find $U$ so that $x_n \in \pi_0 (\bR F(\Aln))$ comes from some $y_n \in \pi_0 (\bR U(\Aln))$, which gives us the bottom two horizontal maps. Using the fact that the Postnikov tower for $A$ is built from square zero extensions, one can show that for large $n$ the morphism $\bR U(A_{\leq n+1}) \rarr \bR U(\Aln) \times^h_{\bR F(\Aln)} \bR F(A_{\leq n+1})$ is surjective, which provides us with the bottom square, hence the dotted maps from $\holim \bR U(\Aln)$, as well as the horizontal dotted map between homotopy limits. The map $\bR U(A) \rarr \holim \bR U(\Aln)$ is an equivalence since $U$ is representable, the top horizontal map we have by functoriality, and the right diagonal map is the natural map, which we aim to show is an equivalence. By an inductive argument, which is the statement of Lemma C.0.10, the top square is homotopy cartesian, which allows us to conclude. Thus it seems everything hinges on the local existence of $U$.\\

This we have by virtue of Corollary 2.2.2.9, which essentially says that the convergence of the tower for $\bR F(A)$ can be studied locally on the small \'etale site of $A$ as observed in \cite{TV4}, which is an analogous statement to having $X \rarr Y$ $\rho + 1$-connected in Goodwillie calculus. Both statements provide us with a notion of neighborhood. Further in the above diagram, one can climb up in the Postnikov tower of $\bR U(A)$ for large values of $n$ by Corollary 2.2.5.3. of \cite{TV4}, the argument being that some fibers are simply connected, which immediately gives us the desired convergence, and this also mirrors the fact that the connectedness of $F(X) \rarr P_nF(X)$ tends to zero as $n$ goes to infinity for the convergence of the Taylor tower of homotopy functors. Both theories use deformation theoretic arguments, the existence of a cotangent complex for stacks, and the fact that the Postnikov tower for $A \in \skAlg$ is built up from iterated square zero extensions, and the use of connectedness arguments in the case of Goodwillie calculus, along with repeated joins from $X$, providing us with $T_n^i F(X)$, the homotopy limit of which is $P_n F(X)$.\\

Thus heuristically an overarching theory would go like so: let $\cC$ be a category with a notion of neighborhood, let $\cD$ be a category with a notion of equivalence, $F: \cC \rarr \cD$ a functor, $X \in \Ob \cC$, $\cdots \rarr X_n \rarr X_{n-1} \rarr \cdots$ some deformation theoretic tower induced by $X$. If $F$ satisfies some deformation theoretic constraint and $X$ lies in a certain neighborhood, one has $F(X) \cong \holim F(X_n)$ in $\Ho(\cD)$. The abstraction of such a formalism we denote by $\thT_0$.\\

Applying this to Goodwillie Calculus, governed by a theory $\bG$, $\cC$ is a category of spaces over some fixed space $Y$, $\cD$ is a category of spaces or spectra. On $\cC$ the notion of neighborhood is provided by connectedness: $X \rarr Y$ is $\rho + 1$-connected. $F: \cC \rarr \cD$ is a homotopy functor. The tower induced by $X$ is constructed using the join $X *_Y U$ and its iterations, giving rise first to $T_nF(X)$, then its powers $T_n^i F(X)$, and finally $P_nF(X)$. The constraint put on $F$ is being $\rho$-analytic, which controls the connectedness of the image of cubical diagrams by $F$, and in particular those of the form $F(X) \rarr T_nF(X)$. That we use connectedness is natural since $F$ is valued in a category of spaces, on which weak homotopy equivalence is the natural notion of equivalence. Thus $\bG$ presents itself as an expansion of $\thT_0$.\\

For geometric stacks in derived Algebraic Geometry, governed by a theory $\bS$, $\cC = \skAlg$, $\cD = \SetD$. However, the functor $F$ is a stack, object of $\kDAfftildeet$, thus the notion of neighborhood on $\cC$ is really given by the fact that we have an \'etale model topology on $\kDAff$. This is illustrated in Corollary 2.2.2.9 of \cite{TV4} which shows that one can work locally on the small \'etale site of $A$. That provides us with a neighborhood of $A$. The tower for $A$ is given by its Postnikov tower in $\skAlg$. It is deformation theoretic in nature since one can show it is built from square zero extensions. The constraint we put on $F$ is that it be an op-excisive $n$-geometric $D^{-}$-stack, which implies in particular that it has an obstruction theory, which we saw is the analog of being $\rho$-analytic for homotopy functors. To have an obstruction theory means in particular that $F$ has a cotangent complex. This is necessary to show that, along with the fact that locally on the small \'etale site of $A$ (which guarantees we have a smooth morphism $U \rarr F$ from some representable stack $U$), the homotopy fiber of $\bR U(A_{\leq n+1}) \rarr \bR U(\Aln) \times^h_{\bR F(\Aln)} \bR F(A_{\leq n+1})$ is simply connected for large values of $n$, we have the isomorphism $\bR F(A) \cong \holim \bR F(\Aln)$. Thus $\bS$ presents itself as an expansion of $\thT_0$ as well.\\

\bigskip
\footnotesize
\noindent
Renaud Gauthier, \textsc{University of Mary, 7500 University Dr., Bismarck, ND 58504, USA} \par \nopagebreak \noindent \textit{e-mail address}: \texttt{rg.mathematics@gmail.com}.

\end{document}